\documentclass[12pt]{article}
\usepackage{amssymb,amsmath,amsthm,mathrsfs,txfonts,color}
\usepackage[numbers,sort&compress]{natbib}
\newtheorem{theorem}{\indent Theorem}[section]
\newtheorem{lemma}{\indent Lemma}[section]

\newtheorem{definition}{\indent Definition}[section]
\newtheorem{remark}{\indent Remark}[section]

\let\oldsection\section
\renewcommand\section{\setcounter{equation}{0}\oldsection}

\allowdisplaybreaks

\begin{document}
\title{A reducing mechanism on wave speed for chemotaxis systems
with degenerate diffusion}
\author{
Shanming Ji$^{a}$,
Zhian Wang$^{b}$,
Tianyuan Xu$^{c,d}$\thanks{Corresponding author, email:tian.yuan.xu@mail.mcgill.ca},
Jingxue Yin$^{c}$,
\\
\\
{ \small \it $^a$School of Mathematics, South China University of Technology}
\\
{ \small \it Guangzhou, Guangdong, 510641, P.~R.~China}
\\
{ \small \it $^b$Department of Applied Mathematics, Hong Kong Polytechnic University}
\\
{ \small \it Hung Hom, Kowloon, Hong Kong}
\\
{ \small \it $^c$School of Mathematical Sciences, South China Normal University}
\\
{ \small \it Guangzhou, Guangdong, 510631, P.~R.~China}
\\
{ \small \it $^d$Department of Mathematics and Statistics, McGill University}
\\
{ \small \it Montreal, Quebec,   H3A 2K6, Canada}
}
\date{}
\maketitle

\noindent {\bf Abstract}. This paper is concerned with traveling wave solutions for a chemotaxis model with degenerate diffusion of porous medium type. We establish the existence of semi-finite traveling waves, including the sharp type and $C^1$ type semi-finite waves. Our results indicate that chemotaxis slows down the wave speed of semi-finite traveling wave, that is, the traveling wave speed for chemotaxis with porous medium (degenerate) diffusion is smaller than that for the porous medium equation without chemotaxis. As we know, this is a new result not shown in the existing literature. The result appears to be a little surprising since chemotaxis is a connective force.  We prove our results by the Schauder's fixed point theorem and estimate the wave speed by a variational approach.

\medskip

\noindent \textbf{Keywords}:
Chemotaxis, Degenerate Diffusion, Semi-finite Waves, Wave Speed.

\noindent \textbf{Mathematics Subject Classification}:  35K65, 35K57, 92C17.

\section{Introduction}
Chemotaxis, the movement of an organism in response to a chemical signal, has been an important mechanism accounting for various biological phenomena/processes, such as aggregation of bacteria, slime mould formation, fish pigmentation, tumor angiogenesis,  blood vessel formation, wound healing  and so on. A general form of chemotaxis models can be encapsulated in the following system
\begin{equation}\label{KSG}
\begin{cases}
\displaystyle
u_t=\nabla\cdot(D(u)\nabla u-\chi u\nabla \phi(v))+f(u,v), \\[1mm]
\tau v_t=\Delta v+h(u,v),
\end{cases}
\end{equation}
where $u(x,t)$ and $v(x,t)$ denote the cell density and chemical (signal) concentration at position $x$ and time $t$, respectively. The function $D(u)$ denotes the cell diffusion rate and  {\normalsize$\phi(v)$} is called the chemotactic sensitivity function describing the signal response mechanism; $f(u,v)$ and $h(u,v)$ denote the birth-death dynamics of the cell and  the chemical, respectively. The constant {\normalsize$\tau \geq 0$} is linked to the diffusion rate of the chemical, and $\chi>0$ is referred to as the chemotactic coefficient  measuring the strength of the chemotactic response.

Depending on the specific forms of $D(u), \phi(v), f(u,v)$ and $h(u,v)$, the model (\ref{KSG}) has generic applications among which the mostly studied case in the literature is that $D(u)$ is constant.  There are two major classes of chemotactic response function: logarithmic response $\phi(v)=\ln v$ and  linear response $\phi(v)=v$.
In the case of logarithmic response, the model \eqref{KSG} with $D(u)=1, f(u,v)=0, \ h(u,v)=-uv^m (m\geq 0)$ was originally proposed in \cite{KS3} to model the wave propagation of bacterial chemotaxis describing the movement of bacteria towards the higher concentration of nutrient consumed by bacteria. The mostly prominent feature of  this type model is the use of logarithmic sensitivity function $\phi(v)=\ln v$, which results in a mathematically unfavorable singularity but is necessary to generate traveling wave solutions (cf. \cite{LuiWang}). This model has inspired a large number of  mathematical works on the study of chemotactic traveling waves (cf. \cite{Choi1, Choi2, Davis, li2009nonlinear, LiLiWang, LiWang2011JDE, Wang-review}) as well as global well-posedness (cf. \cite{Hou1, LiZhao, Martinez} and references therein). In the case of linear response, the prototypical form of \eqref{KSG} was $D(u)=1, f(u,v)=0, h(u,v)=u-v$ originally proposed in \cite{KS1, KS2}, called the minimal chemotaxis model, to model the self-aggregation of {\it Dictyostelium discoideum} in response to cyclic adenosine monophosphate (cAMP) secreted by themselves. The prominent feature of the minimal chemotaxis model is the blowup of solutions depending on the space dimensions (cf.  \cite{Horstman, winkler2013} and references therein).
Collapsing steady states (cf. \cite{NonlinearityPino}) and large mass boundary condensation patterns (cf. \cite{JDEPino})
are established for the steady minimal chemotaxis model. Moreover a number of modifications of the minimal chemotaxis model have been proposed so that the modified models allow pattern formation without blow-up (see a survey paper \cite{hp2009}).

In the above-mentioned works, the diffusion rate $D(u)$ is assumed to be constant which neglects the volume exclusion effect of packed cells in chemotactic process. By taking this effect into account, several types of density-dependent diffusion
have been introduced to  chemotaxis models, among which is the following chemotaxis model with porous medium type diffusion (cf. \cite{hp2009,Kowalczyk})
\begin{equation}\label{eq-dc}
\begin{cases}
\displaystyle
u_t=\Delta u^m-\chi\nabla\cdot(u\nabla v)+u(a-bu),\\[1mm]
\tau v_t=\Delta v-v+u.
\end{cases}
\end{equation}
which, compared to the minimal chemotaxis system, introduces the effect of population pressure ($m>1$) and cellular growth $u(a-bu)$ (if $a,b>0$). The use of degenerate diffusion is very important in mathematical modeling of biological processes. Typically, invasion of migrating populations is characterized by a well-defined boundary. This feature can not be reproduced by linear diffusion which gives rise to the smooth fronts being positive everywhere.
However, the degenerate diffusion raises the possibility of sharp-front waves relevant to the advancing distinct edge of an expending cell population, such as a growing tumour or moving bacterial colonies.
Numerous results have been developed for the chemotaxis system \eqref{eq-dc} in bounded domain with Neumann boundary conditions. When $m=1$, it was first shown in \cite{Tello} that the globally bounded solutions of \eqref{eq-dc} exist in two dimension or higher dimensions if $b>0$ is sufficiently large. Many refined results for $m=1$ have been developed afterwards (e.g. see \cite{Xiang, winkler2010boundedness} and references therein). Without logistic growth (i.e. $a=b=0$), the globally bounded weak solutions of \eqref{eq-dc} with $m>1$  have been obtained in two or higher dimensions (cf. \cite{Blanchet, Ishida, Sugiyama} and references therein) and degenerate stationary solutions have been extensively studied in \cite{Carrillo,Carrillo16,Kim}. In this paper, we are concerned with the wave propagation dynamics of \eqref{eq-dc}.

Without cell growth, it has been shown in \cite{LuiWang} that the chemotactic sensitivity function $\phi(v)$ in \eqref{KSG} must have a singularity like the logarithmic sensitivity used in \cite{KS3} to generate traveling wave solutions.
However, with cell growth (i.e. proliferative cell population),  traveling wave solutions are possible without singular chemotactic sensitivity.
In recent works, Salako and Shen \cite{salako2016spreading,salako2018existence} studied the existence of traveling
waves of parabolic-elliptic chemotaxis system \eqref{eq-dc} (i.e.  $\tau=0$) with linear diffusion $m=1$ in $\mathbb{R}$.
Then, they extended their results to parabolic-parabolic cases (i.e. $\tau=1$) \cite{Shen19}.
We noted that the results of \cite{salako2016spreading,salako2018existence,Shen19} showed that the minimal wave speed is  $2\sqrt{a}$, which is the same as the one for the well-known Fisher-KPP equation (namely the first equation of \eqref{eq-dc} with $\chi=0$ and $m=1$). This indicates that for the linear diffusion ($m=1$), the chemotaxis does not affect the traveling wave speed. In this paper, we shall show that once $m>1$, the chemotaxis will reduce the minimal wave speed. This is a new finding not observed in the literature as far as we know.

In the absence of chemotaxis (i.e. $\chi=0$), the first equation of \eqref{eq-dc} becomes the following (degenerate) Fisher-KPP equation
\begin{equation}\label{FKPP}
u_t=\Delta u^m+u(a-bu),\quad x\in \mathbb R^N.
\end{equation}
Without growth term ($a=b=0$), the equation \eqref{FKPP} is well-known as the porous medium equation if $m>1$, which has significant applications in physics and been extensively studied over many years. For instance,   M. del Pino et al. studied the solvability, extinction profiles,
optimal decay rates \cite{PinoARMA, PinoIndiana,PinoJMPA} for the fast diffusion case ($m<1$) and
 V\'azquez analyzed the exact behaviors of solutions for the case porous medium case ($m>1$) in \cite{VazquezOxford,VazquezOxfordSmooth},
such as asymptotics, propagation properties, and stability.
With growth term ($a,b>0$), \eqref{FKPP} is called the classical Fisher-KPP equation if $m=1$ and degenerate Fisher-KPP equation if $m>1$.
It is widely known that the classical Fisher-KPP equation has a unique (up to a translation) positive traveling wave solution with a minimum wave speed, while in contrast the degenerate Fisher-KPP equation allows for a compactly supported traveling wave with a sharp edge. The feature of this sharp type semi-finite traveling wave is that its profile has distinct boundaries and decays to zero at a finite point in space rather than decaying to zero at far field as for the case of linear diffusion ($m=1$). The degenerate Fisher-KPP equation ($m>1$) has important applications in describing biological processes, such as the spreading of MG63 cells \cite{sengers2007experimental}, cell migration during epidermal wound healing \cite{sherratt1990models}, ``population pressure'' in adaptive bacterial colonies \cite{Ben1994Generic}, and so on.
It was Aronson \cite{Aronson} who first discovered this sharp type semi-finite traveling waves for the degenerate Fisher-KPP equation and found an interesting phenomenon: the traveling wave with the critical wave speed $c_*(m,0)$ is the only one with semi-finite type. The similar behavior has been found for the reaction diffusion equation with degenerate diffusion and convection (cf. \cite{gilding2005fisher}) and degenerate diffusion model for bacterial pattern formation (cf. \cite{satnoianu2001travelling}). Further results are obtained in Audrito and V\'azquez \cite{VazquezJDE17} for a doubly nonlinear diffusion equation and in \cite{HJMY} with a variational characterization of critical wave speed for other type reaction terms. Except important applications in physics and biology as mentioned above, the porous medium diffusion has also been used to describe some dynamics in ecology and animal dispersal (cf. \cite{Mathysen,murray2002mathematical}).

In this paper, we shall consider the chemotaxis system \eqref{eq-dc} with generate diffusion ($m>1$), which can also be regarded as the Fisher-KPP model with chemotactic advection ($\chi>0$).
Though there are  many results on various chemotaxis systems with density-dependent diffusion (cf. \cite{Carrillo,Carrillo16, Kim, Blanchet, Ishida, Sugiyama, newgreen2003chemotactic, xuMBE, xuM2as, xuarxiv19}), and there are no results on the traveling wave solutions. The aim of this paper is to examine the existence of traveling wave solutions of \eqref{eq-dc} with $m>1$ and in particular we shall investigate how the chemotaxis will affect the wave propagating properties compared to the degenerate Fisher-KPP equation \eqref{FKPP} and how the degeneracy of diffusion will make difference compared to the results obtained in \cite{salako2016spreading,salako2018existence,Shen19} for  \eqref{eq-dc} with linear diffusion $m=1$. Roughly speaking, we prove that \eqref{eq-dc} admits a monotonically increasing semi-finite type traveling wave solution in $\mathbb{R}$ with speed $c_*(m,\chi)\in[2c_*(m,0)/3,c_*(m,0))$ where $c_*(m,0)$ is the wave speed of the unique semi-finite traveling wave for degenerate Fisher-KPP equation \eqref{FKPP} as in \cite{Aronson}, see Theorem \ref{th-main}. Clearly the wave speed $c_*(m,\chi)$  is smaller than that for the degenerate Fisher-KPP equation \eqref{FKPP}. This essentially shows that chemotaxis is a factoring slowing down the propagating wave speed, which is a new result not discovered in the literature. Compared to the results obtained in \cite{salako2016spreading,salako2018existence,Shen19} for \eqref{eq-dc} with linear diffusion $m=1$, our results imply that the degenerate diffusion will substantially change the traveling wave profile and wave speed. Moreover we show that the semi-finite traveling wave solutions obtained in Theorem \ref{th-main} is of sharp type if $m\geq 2$ and $C^1$ type if $1<m<2$ (see Theorem \ref{th-sharp}). This shows the effect of degeneracy of diffusion on the regularity of semi-finite traveling wave solutions.  To prove our results, we develop a new framework by using Schauder's fixed point theorem to prove the existence of semi-finite traveling waves of the degenerate chemotaxis model \eqref{eq-dc} in $\mathbb{R}$ for which the maximum principle is inapplicable. Due to the complexity of the model structure, it is very difficult to determine the explicit speed of travelling waves.
In this paper, we first propose a variational characterization to the degenerate chemotaxis model to estimate the wave speed.

The organization of this paper is as follows. In Section 2, we state our main results.
In section 3, we shall prove the existence of semi-finite traveling waves solutions to the system \eqref{eq-dc}
and show the estimation of the wave speed.

\section{Statement of main results}

We consider the existence of semi-finite type traveling waves of the
chemotaxis model with degenerate diffusion \eqref{eq-dc} on real line $\mathbb{R}$. In the sequel, the space is always fixed with $\mathbb{R}$ without emphasis any more.
Let $\xi=x+ct$ with speed $c>0$ be the moving coordinate.
The traveling wave $u(x,t)=\phi(\xi)$, $v(x,t)=\eta(\xi)$
with $\xi=x+ct$ of \eqref{eq-dc} satisfies
\begin{equation}\label{eq-TW}
\begin{cases}
-(\phi^m)_{\xi\xi}+c\phi_\xi+\chi (\eta_\xi\phi)_\xi
=\phi(a-b\phi), \quad \xi\in\mathbb R,\\
-\eta_{\xi\xi}+\tau c\eta_\xi+\eta=\phi, \quad \xi\in\mathbb R,\\
\phi(-\infty)=\eta(-\infty)=0,\quad \phi(+\infty)=\eta(+\infty)=a/b.
\end{cases}
\end{equation}
The second equation in \eqref{eq-TW} is a linear diffusion equation
and can be solved by the following Bessel potential
\begin{equation} \label{eq-bessel}
\displaystyle
\eta(\xi;\phi)=\int_0^\infty\int_{\mathbb R}
\frac{e^{-s}}{\sqrt{4\pi s}}e^{-\frac{|\xi-z|^2}{4s}}\phi(z-\tau cs)dzds
=:B[\phi](\xi).
\end{equation}
Therefore, the profile function $\phi(\xi)$ satisfies
\begin{equation} \label{eq-u2}
-(\phi^m)_{\xi\xi}+c\phi_\xi+\chi ((B[\phi])_\xi\phi)_\xi
=\phi(a-b\phi),
\end{equation}
or equivalently,
\begin{equation} \nonumber
-(\phi^m)_{\xi\xi}+(c+\chi (B[\phi])_\xi)\phi_\xi
=\phi(a-\chi(B[\phi]+\tau c(B[\phi])_\xi)-(b-\chi)\phi),
\end{equation}
which can be regarded as a second order degenerate differential equation
with nonlocal convection caused by the chemotaxis.

One of the important features of the degenerate diffusion equation is
that it generates compactly supported solutions for
compactly supported initial values.
To reveal this peculiar phenomenon,
we are looking for traveling waves of semi-finite type
and for the sake of convenience we present the following definition of
semi-finite and positive traveling waves.

\begin{definition} \label{de-semi}
A profile function $\phi(\xi)$ is said to be a traveling wave solution (abbreviated as TWS hereafter)
of \eqref{eq-u2} (and also of \eqref{eq-TW} together with $\eta=B[\phi]$)
if $\phi\in C_\text{unif}^b(\mathbb R)$, $0\le \phi(\xi)\le a/b$,
$\phi(-\infty)=0$, $\phi(+\infty)=a/b$,
$\phi^m\in W_{\text{loc}}^{1,2}(\mathbb R)$,
$\phi(\xi)$ satisfies \eqref{eq-u2} in the sense of distributions.
The TWS $\phi(\xi)$ is said to be of semi-finite type if
the support of $\phi(\xi)$ is semi-finite, i.e.,
$\text{supp}\,\phi=[\xi_0,+\infty)$ for some $\xi_0\in\mathbb R$,
$\phi(\xi)>0$ for $\xi>\xi_0$.
On the contrary, the TWS $\phi(\xi)$
is said to be of positive type if $\phi(\xi)>0$ for all $\xi\in\mathbb R$.

Furthermore, for the semi-finite TWS $\phi(\xi)$,
if $\phi''\not\in L_\text{loc}^1(\mathbb R)$,
we say that $\phi(\xi)$ is a sharp type semi-finite TWS;
otherwise, if $\phi''\in L_\text{loc}^1(\mathbb R)$,
we say that $\phi(\xi)$ is a $C^1$ type semi-finite TWS.
\end{definition}

To show our main results on the chemotaxis model with degenerate diffusion in this paper,
we need to recall the following two relevant results concerned with
degenerate diffusion Fisher-KPP equation without chemotaxis (i.e. $m>1$ and $\chi=0$)
and chemotaxis model with linear diffusion (i.e. $m=1$ and $\chi>0$).

\begin{theorem}[TWS for Degenerate Fisher-KPP equation \cite{Aronson}] \label{th-dege}
For the degenerate Fisher-KPP equation \eqref{FKPP} with $m>1$, there exists a critical wave speed $c_*=c_*(m,0)>0$
only depending on $m,a,b$, such that: there are no admissible TWS for $0<c<c_*(m,0)$,
while there exists exactly one admissible TWS for all $c\ge c_*(m,0)$. Moreover
the TWS corresponding to $c_*(m,0)$ is semi-finite,
while for $c>c_*(m,0)$ the TWS are positive.
\end{theorem}

As we know it was Aronson \cite{Aronson} who first proved the above Theorem \ref{th-dege}
about the admissible traveling wave speeds.
Those properties indicate an interesting phenomenon of degenerate diffusion equations:
the TWS corresponding to the critical wave speed
is the only one that is of semi-finite type.
Further results are obtained in Audrito and V\'azquez \cite{VazquezJDE17}
for a doubly nonlinear diffusion equation
and in \cite{HJMY} with a variational characterization
of $c_*(m,0)$ for Nicholson's blowflies model.
We can verify that the same variational characterization in \cite{HJMY} holds true
for the Fisher-KPP equation \eqref{FKPP}.
Therefore, the critical wave speed $c_*(m,0)$ for Fisher-KPP equation \eqref{FKPP}
with $m>1$ is determined by
\begin{equation} \label{eq-cstar0}
c_*(m,0)=\sup_{g\in \mathscr{D}}\left\{2\int_0^{a/b}
\sqrt{-ms^{m-1}g(s)g'(s)(as-bs^2)}ds\right\}=:\sup_{g\in \mathscr{D}}J(g),
\end{equation}
where $\mathscr{D}=\{g\in C^1([0,{a/b}]);g({a/b})=0,
\int_0^{a/b} g(s)ds=1,g'(s)<0,\forall s\in(0,{a/b})\}$.
The functional $J(g)$ attains its maximum on $\mathscr{D}$
at some $\tilde g(\cdot)\in\mathscr{D}$
(see Lemma \ref{le-functional})
and the function $\tilde g$ is only dependent on $m$ and $a,b$.
We denote
\begin{equation} \label{eq-sigma}
\sigma:=\tilde g(0)>0,
\end{equation}
which is a positive constant only depending on $m$ and $a,b$
as proved in Lemma \ref{le-functional}.

\begin{theorem}[Chemotaxis model with linear diffusion \cite{Shen19}] \label{th-chemo}
For the chemotaxis model \eqref{eq-dc} with linear diffusion (i.e. $m=1$),
if $0<\chi<b/2$ and $\tau\ge(1-1/a)_+/2$,
then there exists a TWS for every $c\ge 2\sqrt{a}$;
while there is no admissible TWS for every $c<2\sqrt{a}$.
\end{theorem}

We note that for the classical Fisher-KPP equation
\eqref{FKPP} with $m=1$, it is  well-known that the minimal wave speed is  $2\sqrt{a}$.
Therefore, Theorem \ref{th-chemo} shows that
for the chemotaxis model \eqref{eq-dc} with linear diffusion, the chemotaxis has no effect on the minimal wave speed
under the conditions therein.

The main purpose of this paper is to show the existence of semi-finite type TWS for the
chemotaxis model \eqref{eq-u2} with degenerate diffusion (i.e. $m>1$ and $\chi>0$),
and more importantly, the chemotaxis does slow down the wave speed. That is, the speed of semi-finite TWS is smaller than
that without chemotaxis (which is $c_*(m,0)$ according to Theorem \ref{th-dege}).
\medskip

Our first main result is stated as follows:
\begin{theorem}\label{th-main}
For the chemotaxis model \eqref{eq-dc} with degenerate diffusion ($m>1$),
we assume that
$$
0<\chi\le\min\left\{\frac{2b^2c_*(m,0)}{3a^2\sigma},
\frac{2bc_*(m,0)}{3a},\frac{b}{\tau c_*(m,0)+2}
\right\},
$$
where $c_*(m,0)$ is the critical wave speed of the
degenerate Fisher-KPP equation given in Theorem \ref{th-dege}
and $\sigma$ is a positive constant defined in \eqref{eq-sigma}.
Then there exists a wave speed $c_*(m,\chi)\in[2c_*(m,0)/3,c_*(m,0))$
such that \eqref{eq-u2} admits a monotonically increasing
semi-finite type TWS $\phi(\xi)$ with speed $c_*(m,\chi)$.
\end{theorem}

\begin{remark}
{\em The known results from Aronson \cite{Aronson}, Audrito and V\'azquez \cite{VazquezJDE17}
and references therein for the degenerate diffusion equation without chemotaxis (i.e. \eqref{eq-dc} with $\chi=0$).
show that the semi-finite TWS is unique (up to translation and reflection)
and hence the semi-finite wave speed $c_*(m,0)$ is unique,
which is also the minimal admissible wave speed.
Here, we show that the semi-finite wave speed $c_*(m,\chi)$
for the degenerate diffusion equation with chemotaxis
is smaller than $c_*(m,0)$.
That is, the chemotaxis slows down the speed of the semi-finite type TWS.
We conjure that $c_*(m,\chi)$ is unique and also is the minimal admissible wave speed
for the chemotaxis model with degenerate diffusion.
}
\end{remark}

\begin{remark}{\em
Since $\phi(\xi)$ satisfies \eqref{eq-u2},
which is non-degenerate within $\text{supp}\,\phi$,
it should be noted that the possible singularity of $\phi''$
only occurs at $\{\xi_0\}$, i.e.,
there may be a Radon measure supported in $\{\xi_0\}$.
Therefore, for the sharp type semi-finite TWS $\phi(\xi)$,
we have $\phi(\xi)\not\in C^1(\mathbb R)$;
while for the $C^1$ type semi-finite TWS,
we have $\phi(\xi)\in C^1(\mathbb R)$.
}
\end{remark}

The semi-finite traveling wave solutions are classified into
sharp type and $C^1$ type according to the degeneracy index $m$, as stated in the following theorem.

\begin{theorem} \label{th-sharp}
Assume that the conditions in Theorem \ref{th-main} hold.
If $m\ge2$, then the semi-finite TWS obtained in Theorem \ref{th-main} is of sharp type;
while if $1<m<2$, then the semi-finite TWS is of $C^1$ type.
\end{theorem}

\begin{remark}{\em
There is no semi-finite TWS for the linear diffusion case $m=1$ and hence
roughly speaking the degeneracy is enforced by $m>1$.
The regularity of the case $m\ge2$ is weaker than that of $1<m<2$
and sharp type TWS arises for $m\ge2$ while
there is no sharp type TWS for $1<m<2$.
For the case $1<m<2$, the semi-finite TWS remains $C^1$ regularity
and even higher regularity  if $1<m<1+\varepsilon$ with  small $\varepsilon>0$, but none of them are analytic.
}
\end{remark}

For the sake of convenience, we shift $\xi_0$ to $0$
as the traveling wave equation \eqref{eq-TW} is autonomous,
and we may also write the variable $\xi$ as $t$ in following sections.

\section{Proof of main results}

The existence of semi-finite traveling wave solutions
of \eqref{eq-u2} is proved by the
Schauder's fixed point theorem on an appropriate profile set $\Phi$
constructed with upper and lower profiles $\overline \phi$ and $\underline \phi$
to be specified later.
That is,
\begin{align} \nonumber
\Phi=&\{\phi\in C_\text{unif}^\text{b}(\mathbb R);
0\le\phi(t)\le a/b,
\phi(t)=0,~t\le0, ~\phi(t)>0,~t>0,
\\ \label{eq-Phi}
&\phi(+\infty)=a/b, ~
\underline \phi(t)\le\phi(t)\le\overline\phi(t),~
\phi(t) \text{~is monotonically increasing on~}\mathbb R
\}.
\end{align}
For any $\hat \phi\in\Phi$,
let as in \eqref{eq-bessel}
\begin{equation} \label{eq-eta}
\hat\eta(t)=B[\hat\phi](t)
\end{equation}
be the bounded solution of the second equation in \eqref{eq-TW}
corresponding to $\hat\phi(t)$,
and we solve the following degenerate equation
\begin{equation} \label{eq-aux}
\begin{cases}
-(\phi^m)''+c\phi'+\chi (\hat\eta'\phi)'+b\phi^2=a\phi, \quad t\in\mathbb R,\\
\phi(t)=0, ~t\le0, \quad \phi(t)>0, ~t>0, \quad
\phi(+\infty)=a/b.
\end{cases}
\end{equation}
The above problem will be first solved on the half real line $[0,+\infty)$
and then we extend the solution to $(-\infty,0)$ as
$\phi(t)=0$ for $t<0$.
We note that the extended function $\phi(t)$ satisfies
the degenerate diffusion equation \eqref{eq-aux} on the whole line
in the sense of distributions
provided that $(\phi^m)''\in L_\text{loc}^1(\mathbb R)$.

Since $\hat\eta(t)$ in \eqref{eq-eta} solves the linear differential equation
\eqref{eq-TW}$_2$ corresponding to $\hat\phi(t)$,
we can rewrite \eqref{eq-aux} as
$$
-(\phi^m)''+(c+\chi\hat\eta')\phi'+b\phi^2=(a-\chi\hat\eta'')\phi.
$$
Now, for a general case, we consider the solvability of the following
auxiliary degenerate differential problem
\begin{equation} \label{eq-aux-g}
\begin{cases}
-(\phi^m)''+\lambda(t)\phi'+b\phi^2=\mu(t)\phi, \quad t\in\mathbb R,\\
\phi(t)=0, ~t\le0, \quad \phi(t)>0, ~t>0, \quad
\phi(+\infty)=a/b,
\end{cases}
\end{equation}
with $\lambda(t)\in C^1(\overline{\mathbb R})$, $\mu(t)\in C(\overline{\mathbb R})$
and $\mu(+\infty)=a$.

It should be noted that the differential problem
\eqref{eq-aux-g} (or \eqref{eq-aux}) with singularity at $\{0\}$
is not always solvable.
We just mention that for the case $\chi=0$, or equivalently $\lambda(t)\equiv c$ and
$\mu(t)\equiv a$, the only solvable situation is that
$c=c_*(m,0)$ and $\phi(t)$ is the semi-finite TW
as in Theorem \ref{th-dege}.
Therefore, we solve \eqref{eq-aux-g} locally and modify \eqref{eq-aux-g} as follows
\begin{equation} \label{eq-aux-m}
\begin{cases}
-(\phi^m)''+\lambda(t)\phi'+b\phi^2=\mu(t)\phi, \quad t\in(-\infty,t^*),\\
\phi(t)=0, ~t\le0, \quad \phi(t)>0 \ \text{and} \
\phi(t) \text{~is strictly increasing on~} (0,t^*),
\end{cases}
\end{equation}
where $t^*\in(0,+\infty]$ and $(-\infty,t^*)$ is the maximal solvable interval.

Here are some simple but useful estimates on the Bessel potential
$\hat\eta(t)$.

\begin{lemma} \label{le-besell}
For $0\le\hat\phi(t)\le a/b$, and $\hat\phi(t)\not\equiv0,\not\equiv a/b$,
$\hat\phi(+\infty)=a/b$,
the Bessel potential
$\hat\eta(t)=B[\hat\phi](t)$ in \eqref{eq-eta} satisfies
$$
0<\hat\eta(t)<\frac{a}{b}, \quad |\hat\eta'(t)|\le \frac{a}{2b},
\quad |\hat\eta''(t)|\le (\frac{\tau c}{2}+1)\frac{a}{b},
$$
and $\hat\eta(+\infty)=a/b$, $\hat\eta'(+\infty)=\hat\eta''(+\infty)=0$.
If further we assume that $\hat\phi(t)$ is monotonically increasing,
then $\hat\eta'(t)>0$ and $\hat\eta''(t)\ge -a/b$.
\end{lemma}
{\it\bfseries Proof.}
The  first conclusions can be proved with the fundamental calculus,
see for example \cite{salako2018existence}.
We note that
$$
|\hat\eta''(t)|
=|\tau c\hat\eta'+\hat\eta-\hat\phi|
\le \tau c|\hat\eta'|+|\hat\eta-\hat\phi|
\le (\frac{\tau c}{2}+1)\frac{a}{b},
$$
and $\hat\eta''(t)=\tau c\hat\eta'+\hat\eta-\hat\phi>-\hat\phi\ge -a/b$. This completes the proof.
$\hfill\Box$
\medskip

Now, we  solve the degenerate differential problem
\eqref{eq-aux-m} locally step by step.

\begin{lemma} \label{le-solv}
Assume that $0<\lambda(t)\in C^1(\overline{\mathbb R})$,
$\mu(t)\in C(\overline{\mathbb R})$, $\mu(+\infty)=a$.
Then \eqref{eq-aux-m} admits a local semi-finite type solution $\phi(t)$ such that
$$
\phi(t)=K_1t_+^\frac{1}{m-1}+K_2t_+^\frac{m}{m-1}+o(t_+^\frac{m}{m-1}),
\quad t\to 0,
$$
where $t_+=\max\{t,0\}$, and
$$
K_1=\Big(\frac{m-1}{m}\lambda(0)\Big)^\frac{1}{m-1},
\quad
K_2=-\frac{\mu(0)-\frac{1}{m-1}\lambda'(0)}{2m\lambda(0)}K_1.
$$
For the case $\lambda(t)=c+\chi\hat\eta'(t)>0$ and $\mu(t)=a-\chi\hat\eta''(t)$,
we have
$$
K_1=\Big(\frac{m-1}{m}(c+\chi\hat\eta'(0))\Big)^\frac{1}{m-1},
\quad
K_2=-\frac{a-\frac{m}{m-1}\chi\hat\eta''(0)}{2m(c+\chi\hat\eta'(0))}K_1.
$$
\end{lemma}
{\it\bfseries Proof.}
We first show the existence of the local solution to
the singular differential problem \eqref{eq-aux-m}.
Since $\phi(t)\equiv0$ for $t\le0$, the solution $\phi(t)$ must satisfy
$\phi(0^+)=0$ and $(\phi^m)'(0^+)=0$.
Otherwise, $-(\phi^m)''$ is the summation of a Radon measure
and a locally integrable function,
\eqref{eq-aux-m} can not hold in the sense of distributions.
On the other hand, if $\phi(0^+)=0$ and $(\phi^m)'(0^+)=0$,
then $-(\phi^m)''$ is locally integrable
and $\phi(t)$ only need to satisfy \eqref{eq-aux-m} locally on the interval $(0,t^*)$.
Now, let $\zeta(t)=\phi^m(t)$ and we rewrite \eqref{eq-aux-m} into
\begin{equation} \label{eq-zaux}
\begin{cases}
-\zeta''(t)+\lambda(t)(\zeta^\frac{1}{m})'(t)
=\mu(t)\zeta^\frac{1}{m}(t)-b\zeta^\frac{2}{m}(t), \quad t\in(0,t^*),\\
\zeta(0)=0,~\zeta'(0)=0, ~ \zeta(t)>0, ~t\in(0,t^*), ~
\zeta(t) \text{~is strictly increasing on~} (0,t^*).
\end{cases}
\end{equation}
Integrating \eqref{eq-zaux} over $(0,t)$ implies
\begin{equation} \label{eq-zaux-i}
\begin{cases}
\displaystyle
-\zeta'(t)+\lambda(t)\zeta^\frac{1}{m}(t)
=\int_0^t\Big(
\mu(s)\zeta^\frac{1}{m}(s)-b\zeta^\frac{2}{m}(s)
+\lambda'(s)\zeta^\frac{1}{m}(s)
\Big)ds, \quad t\in(0,t^*),\\
\zeta(0)=0,~ \zeta(t)>0,~ \zeta'(t)\ge0,~t\in(0,t^*).
\end{cases}
\end{equation}
Since $\zeta^\frac{1}{m}$ is not Lipschitz continuous with respect to $\zeta$,
the unique solvability of the equation $\zeta'(t)=\lambda(t)\zeta^\frac{1}{m}(t)$
with the initial condition $\zeta(0)=0$ is not true in general.
In fact, if the condition ``$\zeta(t)>0$'' is removed,
then \eqref{eq-zaux-i} admits infinitely many solutions such that:
$\zeta_0(t)\equiv0$ is the smallest one,
and there exists a unique maximal solution $\zeta_1(t)$
such that $\zeta_1(t)>0$ for $t\in(0,t^*)$ ($\zeta_1$ is the one that we are looking for).
Moreover, $\zeta_1$ is the limiting function of the approximating problem
of \eqref{eq-zaux-i} with ``$\zeta(0)=0$'' being replaced by ``$\zeta(0)=\varepsilon$''
for $\varepsilon>0$ and tending to zero.

Using a formal expansion of $\phi(t)$ near $0$, we can derive the values of those coefficients.
Here we omit the computations.
For the leading term $\phi_0(t):=K_1t_+^\frac{1}{m-1}$, we note that
$\phi_0\not\in C^1(I_0)$ for some neighborhood $I_0$ of $0$ and $m\ge2$.
However, $\phi_0'(t)=\frac{1}{m-1}K_1t_+^\frac{2-m}{m-1}$,
and $[\phi_0^m(t)]''=K_1^m\frac{m}{(m-1)^2}t_+^\frac{2-m}{m-1}$
in the sense of distributions.
This type of $\phi(t)$ is similar to the Barenblatt solution
of the porous medium equation.
$\hfill\Box$

\begin{lemma} \label{le-bounds}
Suppose that the conditions in Lemma \ref{le-solv} hold
and $\frac{1}{3}c_*(m,0)\le \lambda(t)\le \frac{4}{3}c_*(m,0)$,
where $c_*(m,0)$ is the minimal wave speed in Theorem \ref{th-dege}.
Let $(-\infty,t^*)$ be the maximal solvable interval of \eqref{eq-aux-m}.
Then there exist constants $C_1,C_2>0$ depending on $\|\lambda\|_{C^1}$ and $\|\mu\|_{C}$,
such that
$$
K_1t_+^\frac{1}{m-1}-C_1t_+^\frac{m}{m-1}\le
\phi(t)\le K_1t_+^\frac{1}{m-1}+C_2t_+^\frac{m}{m-1},
\quad t<t_1,
$$
where $t_1=t^*$ if $\phi(t^*)\le a/b$ and $t_1=\phi^{-1}(a/b)$ if $\phi(t^*)>a/b$.
\end{lemma}
{\it\bfseries Proof.}
The leading term $\phi_0(t):=K_1t_+^\frac{1}{m-1}$ comes from Lemma \ref{le-solv}
and the tail term follows from the continuous dependence of the coefficients
since $\lambda(t)$ is uniformly bounded in $[\frac{1}{3}c_*(m,0),\frac{4}{3}c_*(m,0)]$
and $\phi(t)$ is strictly increasing on $(0,t^*)$.
$\hfill\Box$

\begin{lemma} \label{le-bounds-psi}
Suppose that the conditions in Lemma \ref{le-solv} are valid
and that $\frac{1}{3}c_*(m,0)\le \lambda(t)\le \frac{4}{3}c_*(m,0)$,
where $c_*(m,0)$ is the minimal wave speed in Theorem \ref{th-dege}.
Let $(-\infty,t^*)$ be the maximal solvable interval of \eqref{eq-aux-m}
and
$$
\psi(t):=(\phi^m(t))'=m\phi^{m-1}(t)\phi'(t), \quad t\in(-\infty,t^*).
$$
We change the variable $t\in(0,t^*)$ as an inverse function of $\phi=\phi(t)$ since
$\phi(t)$ is strictly increasing on $(0,t^*)$
and set $\tilde\psi(\phi)=\psi(t)$.
Then there exist constants $C_3,C_4>0$ depending on $\|\lambda\|_{C^1}$ and $\|\mu\|_{C}$,
such that
$$
\lambda(0)\phi-C_3\phi^2\le
\tilde\psi(\phi)\le \lambda(0)\phi+C_4\phi^2,
\quad \phi\in(0,\phi_1),
$$
where $\phi_1=\frac{a}{2b}$ if $\phi(t^*)\ge\frac{a}{2b}$
and $\phi_1=\phi(t^*)$ if $\phi(t^*)<\frac{a}{2b}$.
\end{lemma}
{\it\bfseries Proof.}
In the interval $(0,t^*)$, according to the relation $\phi=\phi(t)$,
we write $t=\tilde t(\phi)$ and $\tilde\lambda(\phi)=\lambda(\tilde t(\phi))$,
$\tilde\mu=\mu(\tilde t(\phi))$.
The second order degenerate differential equation \eqref{eq-aux-m}
is converted to the following non-autonomous system
\begin{equation} \label{eq-ztilde}
\begin{cases}
\displaystyle
\frac{d\phi}{dt}=\frac{\psi(t)}{m\phi^{m-1}(t)},\\
\displaystyle
\frac{d\psi}{dt}=\lambda(t)\frac{\psi(t)}{m\phi^{m-1}(t)}
+b\phi^2(t)-\mu(t)\phi(t).
\end{cases}
\end{equation}
Then using $\phi$ as the variable, we have
\begin{equation} \label{eq-ztildepsi}
\frac{d\tilde\psi}{d\phi}
=\tilde\lambda(\phi)+\frac{m\phi^{m-1}(b\phi^2-\tilde\mu(\phi)\phi)}{\tilde\psi},
\end{equation}
which is regarded as a generalized phase plane,
see \cite{HJMY} and \cite{JDE18XU} for example.
For the constant case such that $\lambda(t)\equiv c$ and $\mu(t)\equiv a$,
the singular differential system \eqref{eq-ztilde}
admit an unique local semi-finite type solution
such that $\phi(t)=0$ for $t\le0$ and $\phi(t)>0$ for $t\in(0,t^*)$
and infinitely many local positive type solutions
(such that $\phi(t)>0$ for $t\in(-\infty,t^*)$).
Furthermore, the semi-finite type solution
corresponds to a trajectory $\tilde\psi(\phi)$ in the phase plane \eqref{eq-ztildepsi}
such that $\tilde\psi(\phi)\sim c\phi$ as $\phi\to 0^+$;
while the positive type solutions corresponds to infinitely many
trajectories $\tilde\psi(\phi)$ in the phase plane \eqref{eq-ztildepsi}
such that $\tilde\psi(\phi)\sim \frac{am}{c}\phi^m$ as $\phi\to 0^+$.
Detailed analysis can be found in Audrito and V\'azquez \cite{VazquezJDE17}
and in Huang et al. \cite{HJMY} for generalized phase plane analysis.
The trajectory $\tilde\psi(\phi)$ corresponding to a semi-finite type solution
that we are looking for is the maximal solution of
\eqref{eq-ztildepsi} with the initial condition $\tilde\psi(0)=0$,
which can be obtained by solving \eqref{eq-ztildepsi} with
$\tilde\psi_\varepsilon(0)=\varepsilon$ for $\varepsilon>0$
and then sending $\varepsilon\to0^+$.
Similar process is applicable for a general case of \eqref{eq-ztildepsi} with
$\lambda(t)\in C^1$ and $\lambda(t)$ being uniformly bounded above and below.
Therefore,
$$
\tilde\psi(\phi)\sim \lambda(0)\phi, \quad \phi\to0^+.
$$
The uniform estimates on the approximate solutions
$\tilde\psi_\varepsilon(\phi)$ for \eqref{eq-ztildepsi}
with the initial condition $\tilde\psi_\varepsilon(0)=\varepsilon$
completes the proof.
$\hfill\Box$

We need to formulate a range of the wave speed $c_*(m,\chi)$
corresponding to the semi-finite type TWS of \eqref{eq-u2}.
Since our framework is based on the fixed point theorem,
for any given $\hat\phi\in \Phi$,
we must determine the range of the minimal speed $c^*(m,\chi,\hat\phi)$
such that the local semi-finite type solution of \eqref{eq-aux-m}
grows up to $a/b$ as we want.
In the following, we first present an abstract definition of the speed $c^*(m,\chi,\hat\phi)$
and then give upper and lower bound estimates by
using a variational characterization inspired by
\cite{Benguria} for linear diffusion and \cite{HJMY} for degenerate diffusion.

\begin{lemma} \label{le-cstar}
For any $0\le\hat\phi(t)\le a/b$ and $\hat\phi(t)$ monotonically increasing,
let $\hat\eta(t)=B[\hat\phi](t)$ as in \eqref{eq-eta},
$\lambda(t)=c+\chi\hat\eta'(t)$, $\mu(t)=a-\chi\hat\eta''(t)$
and $\phi_c(t)$ be the local semi-finite type solution of \eqref{eq-aux-m}
solved in Lemma \ref{le-solv}.
Then there exists a constant $\overline c>0$
depending on the upper bound of $\chi a/b$, such that
$\phi_c(t)$ strictly increases up to $a/b$ in finite time for $c\ge\overline c$.
Define
$$
c^*(m,\chi,\hat\phi)=\inf\{c>0; \text{~the local semi-finite solution ~}\phi_c(t^*)
\ge {a}/{b}\},
$$
which is well-defined and $c^*(m,\chi,\hat\phi)\in[0,\overline c]$.
\end{lemma}
{\it\bfseries Proof.}
According to Lemma \ref{le-besell}, we see that
$|\chi\hat\eta'(t)|\le \frac{\chi a}{2b}$,
$|\chi\hat\eta''(t)|\le \frac{\chi(\tau c/2+1)a}{b}$,
and $\chi\hat\eta''(t)\ge-\chi a/b$.
Therefore, $\lambda(t)\ge c-\frac{\chi a}{2b}$
and $\mu(t)\le a+\chi a/b$.
We make change of variables as in the proof of Lemma \ref{le-bounds-psi}
and according to \eqref{eq-ztildepsi} we have
\begin{equation} \label{eq-zoverlinec}
\frac{d\tilde\psi}{d\phi}\ge c-\frac{\chi a}{2b}
+\frac{m\phi^{m-1}(b\phi^2-(a+\chi\frac{a}{b})\phi)}{\tilde\psi}.
\end{equation}
Let $\tilde\psi_0(\phi)$ be the maximal solution of the differential problem
of \eqref{eq-zoverlinec} with ``$\geq$'' replaced by ``$=$'' and
the initial condition $\tilde\psi_0(0)=0$.
And the semi-finite type solution corresponding to $\tilde\psi_0(\phi)$
is denoted by $\phi_c^0(t)$.
Similar to the phase plane analysis in Audrito and V\'azquez \cite{VazquezJDE17}
and in Huang et al. \cite{HJMY}, we can show that
$\phi_c^0(t)$ grows up to $a/b$ for all $c\ge \overline c$ with some $\overline c>0$
depending on the upper bound of $\chi a/b$.
The comparison between $\tilde\psi(\phi)$ and $\tilde\psi_0(\phi)$
yields that $\phi_c(t)\ge \phi_c^0(t)$
and $c^*(m,\chi,\hat\phi)\le \overline c$ is well-defined.
The proof is completed.
$\hfill\Box$

Lemma \ref{le-cstar} is insufficient to provide
a lower bound of $c^*(m,\chi,\hat\phi)$.
Here we recall the following variational characterization of the minimal wave speed
in \cite{Benguria} for linear diffusion and \cite{HJMY} for degenerate diffusion.

\begin{lemma}[Variational Characterization of Wave Speed, \cite{Benguria,HJMY}]
\label{le-var}
For the diffusion equation $u_t=(u^m)_{xx}+f(u)$
with $f\in C^1([0,u^+])$, $f(0)=f(u^+)=0$, $f'(0)>0$, and $f(u)>0$ for $u\in(0,u^+)$,
where $u^+>0$ is a positive equilibrium.
Then (i) for the linear diffusion case $m=1$, the minimal admissible wave speed
$$c_*(1,0)=\max\left\{2\sqrt{f'(0)},
\sup_{g\in D}\left\{\frac{2\int_0^{u^+}\sqrt{-f(u)g(u)g'(u)du}}{\int_0^{u^+}g(u)du}
\right\}\right\};$$
(ii) for the degenerate diffusion case $m>1$, the minimal admissible wave speed
$$c_*(m,0)=\sup_{g\in D}
\left\{\frac{2\int_0^{u^+}\sqrt{-mu^{m-1}f(u)g(u)g'(u)du}}{\int_0^{u^+}g(u)du}
\right\},$$
where $D=\{g\in C^1([0,u^+]);g(u^+)=0,g'(u)<0,u\in(0,u^+)\}$.
\end{lemma}

\begin{lemma} \label{le-functional}
The functional
$$
J(g):=2\int_0^{a/b}\sqrt{-mu^{m-1}(au-bu^2)g(u)g'(u)}du
$$
on $\mathscr{D}:=\{g\in C^1([0,a/b]);g(a/b)=0,g'(u)<0,u\in(0,a/b),
\int_0^{a/b}g(u)du=1\}$
attains its maximum at some function $\tilde g\in \mathscr{D}$.
Moreover, $\tilde g(\cdot)$ is only dependent on $m,a,b$ and then $\sigma:=\tilde g(0)$
is a positive constant only depending on $m,a,b$.
\end{lemma}
{\it\bfseries Proof.}
The proof of this lemma could be borrowed from \cite{Benguria,HJMY},
but we propose an alternative proof that provides an estimates on the dependence
of $\tilde g(\cdot)$ on $m$.
Let $g\in D$ be an extreme function of
$$
I(g):=\left\{\frac{2\int_0^{a/b}\sqrt{-ms^{m-1}g(s)g'(s)(as-bs^2)ds}}{\int_0^{a/b}g(s)ds}
\right\},
$$
on the function set
$D:=\{g\in C^1([0,a/b]);g(a/b)=0,g'(s)<0,s\in(0,a/b)\}$.
The variational characterization Lemma \ref{le-var} shows that
$$c_*(m,0)=\sup_{g\in D}I(g)=\sup_{g\in\mathscr{D}}J(g).$$
For any $h\in C_0^1((0,a/b))$, there exists $\varepsilon_0>0$
such that $g+\varepsilon h\in D$ for all $\varepsilon\in(0,\varepsilon_0)$.
We may assume that $a=b=1$ without loss of generality, otherwise we can make change of variables
and similar result holds.
Then
$$
\lim_{\varepsilon\to0}\frac{I(g+\varepsilon h)-I(g)}{\varepsilon}=0.
$$
That is,
\begin{align*}
&\int_0^1\frac{-ms^m(1-s)(g'(s)h(s)+g(s)h'(s))}{\sqrt{-ms^m(1-s)g'(s)g(s)}}ds
\cdot\int_0^1g(s)ds
\\
=&2\int_0^1\sqrt{-ms^m(1-s)g'(s)g(s)}ds\cdot\int_0^1h(s)ds.
\end{align*}
Define $f(s)=ms^m(1-s)$ and
$$\theta(s):=-\frac{g'(s)}{g(s)}>0, \quad s\in(0,1).$$
Then integrating by part as $h\in C_0^1((0,1))$ implies that
$$
\int_0^1\left(\Big(\frac{\sqrt{f}}{\sqrt{\theta}}\Big)'h(s)+\sqrt{f}\sqrt{\theta}h(s)\right)ds
=I(g)\int_0^1h(s)ds.
$$
Noticing that $h\in C_0^1((0,1))$ is arbitrary, we see that
$$
\Big(\frac{\sqrt{f}}{\sqrt{\theta}}\Big)'+\sqrt{f}\sqrt{\theta}=I(g)
$$
and then
\begin{align*}
\Big(\frac{{f}}{{\theta}}\Big)'=
2\Big(\frac{\sqrt{f}}{\sqrt{\theta}}\Big)\Big(\frac{\sqrt{f}}{\sqrt{\theta}}\Big)'
=2I(g)\frac{\sqrt{f}}{\sqrt{\theta}}-2f(s),
\quad s\in(0,1).
\end{align*}

We further define
$$
\rho(s):=\frac{f(s)}{\theta(s)}>0, \quad s\in(0,1).
$$
Then $\rho(s)$ is a solution of the following singular ODE
\begin{equation} \label{eq-zrho}
\rho'(s)=2I(g)\rho^\frac{1}{2}(s)-2f(s),
\quad s\in(0,1).
\end{equation}
If $\rho(1)>0$, then $\theta(1^-)=\frac{f(1^-)}{\rho(1)}=0$
and $|g'(s)|\le g(s)$ for $s\in(1-\delta,1)$ with some $\delta\in(0,1)$.
It follows that $g(s)\equiv0$ since $g(1)=0$, which contradicts $g\in D$.
The singular ODE \eqref{eq-zrho} with $\rho(1)=0$ admits a unique solution
such that $\rho(s)>0$ for $s\in(0,1)$.
In fact, let $\tilde s=1-s$ and $\tilde\rho(\tilde s)=\rho(1-s)$, then
\eqref{eq-zrho} is converted into
\begin{equation} \label{eq-zrho2}
\tilde\rho'(\tilde s)=2m\tilde s(1-\tilde s)^m-2I(g)\tilde\rho^\frac{1}{2}(\tilde s),
\quad \tilde s\in(0,1),
\end{equation}
with $\tilde\rho(0)=0$.
Asymptotic analysis shows that
$\tilde\rho(\tilde s)=A\tilde s^2+o(\tilde s^2)$ as $\tilde s\to 0^+$,
where $A=(\frac{\sqrt{(I(g))^2+4m}-I(g)}{2})^2<m$ and then
$$-\frac{g'}{g}=\frac{m}{A(1-s)+o(|1-s|^2)}, \quad s\to 1^-,$$
which means that $g(s)=C(1-s)^\frac{m}{A}+o(|1-s|^\frac{m}{A})$ for some constant $C>0$.
Furthermore, we have $\tilde \rho(\tilde s)>0$ for $\tilde s\in(0,1)$ and
$\tilde\rho(\tilde s)\ge ({m\tilde s(1-\tilde s)^m}/{I(g)})^2$
for $\tilde s\in(\tilde s_0,1)$
with some $\tilde s_0\in(0,1)$ such that $\tilde \rho$ attains its maximum at $\tilde s_0$
since according to \eqref{eq-zrho2},
$\tilde \rho'(\tilde s)>0$ if $I(g)\tilde \rho^\frac{1}{2}(\tilde s)<m\tilde s(1-\tilde s)^m$,
and $\tilde \rho'(\tilde s)<0$ if $I(g)\tilde \rho^\frac{1}{2}(\tilde s)>m\tilde s(1-\tilde s)^m$.

On the other hand, the solution of
the singular ODE \eqref{eq-zrho} with $\rho(0)=0$ is not unique and
we denote the maximal solution by $\rho_0(s)$
similar to the proof in Lemma \ref{le-solv}.
Asymptotic analysis near $s=0$ shows that
$\rho_0(s)=(I(g))^2s^2+o(s^2)$ as $s\to 0^+$ and
there exists a $s_0\in(0,1)$ depending on $m$ such that
$\rho_0(s)\ge(I(g))^2s^2/2$ for $s\in(0,s_0)$.
The comparison principle shows that $\rho(s)\ge\rho_0(s)$ for $s\in(0,1)$ and then
$$
-\frac{g'(s)}{g(s)}=\theta(s)=\frac{f(s)}{\rho(s)}\le \frac{2ms^m(1-s)}{(I(g))^2s^2},
\quad s\in(0,s_0),
$$
which implies that (noticing that $g(\cdot)$ is strictly decreasing on $(0,1)$)
$$
0<\ln\frac{g(0)}{g(s)}\le\frac{2m}{(I(g))^2}\int_0^st^{m-2}(1-t)dt
\le \frac{2m}{(I(g))^2(m-1)}s^{m-1},
\quad s\in(0,s_0),
$$
since $m>1$.
Now we conclude that $g(s)\ge g(0)e^{-Bs^{m-1}}$ for $s\in(0,s_0)$
with the positive constant $B:=2m/((I(g))^2(m-1))$ and then we normalize $g(s)$ as
$\tilde g(s):=g(s)/\int_0^1g(s)ds$ such that
$$
1=\int_0^1\tilde g(s)ds\ge\int_0^{s_0}\tilde g(0)e^{-Bs^{m-1}}ds
\ge e^{-B}\tilde g(0)s_0.
$$
Therefore, $\tilde g(0)\le e^B/s_0$, which is a constant depending on $m$.
$\hfill\Box$

Inspired by the above variational characterization in
\cite{Benguria} for the linear diffusion and \cite{HJMY} for degenerate diffusion.
We formulate the following variational characterization formula for
the degenerate diffusion equation with chemotaxis.

\begin{lemma} \label{le-varchemo}
For any $0\le\hat\phi(t)\le a/b$ and $\hat\phi(t)$ monotonically increasing,
let $\hat\eta(t)=B[\hat\phi](t)$ as in \eqref{eq-eta},
$c=c^*(m,\chi,\hat\phi)$ as in Lemma \ref{le-cstar},
$\lambda(t)=c+\chi\hat\eta'(t)$, $\mu(t)=a-\chi\hat\eta''(t)$
and $\phi_c(t)$ be the local semi-finite type solution of \eqref{eq-aux-m}
solved in Lemma \ref{le-solv}.
If $\chi\sigma\frac{a^2}{2b^2}\le \frac{1}{3}c_*(m,0)$,
then $t^*=+\infty$ and
$$c^*(m,\chi,\hat\phi)\in\big[\frac{2}{3}c_*(m,0),c_*(m,0)\big),$$
where $\sigma=\tilde g(0)$ is the constant in Lemma \ref{le-functional}.
\end{lemma}
{\it\bfseries Proof.}
For $c=c^*(m,\chi,\hat\phi)$, the locally continuously dependence of $\phi_c(t)$ on the
parameter $c$ implies that
(i) $t^*=+\infty$, $\phi_c(t)$ is increasing on $(0,+\infty)$
and $\phi_c(+\infty)=a/b$; or
(ii) $t^*<+\infty$, $\phi_c(t)$ is increasing on $(0,t^*)$
and $\phi_c(t^*)=a/b$, $\phi_c'(t^*)=0$; or
(iii) $t^*<+\infty$, $\phi_c(t)$ is increasing on $(0,t^*)$
and $\phi_c(t^*)>a/b$, there exists a $t_0\in(0,t^*)$ such that
$\phi_c'(t_0)=0$ and $\phi_c(t_0)<a/b$.
The choice of $t_0$ in case (iii) may not be unique,
we choose $t_0$ to be the smallest one if necessary,
since $\phi_c(t)$ is strictly increasing on $(0,t^*)$
those points are isolated.
In cases (i) and (ii), $\phi_c(t^*)=a/b$ and $\phi_c'(t^*)=0$.

Similar to the proof of Lemma \ref{le-bounds-psi}, we make change of variables
such that for $t\in(0,t^*)$ and $\phi\in(0,a/b)$, where
they are related by $t=\tilde t(\phi)$ by $\phi=\phi_c(t)$.
The second order degenerate differential equation \eqref{eq-aux-m}
is converted to the following non-autonomous system
\begin{equation} \label{eq-ztilde-chi}
\begin{cases}
\displaystyle
\frac{d\phi}{dt}=\frac{\psi(t)}{m\phi^{m-1}(t)},\\
\displaystyle
\frac{d\psi}{dt}=c\frac{\psi(t)}{m\phi^{m-1}(t)}
+\chi(\hat\eta'_t(t)\phi(t))'_t
+b\phi^2(t)-a\phi(t),
\end{cases}
\end{equation}
where $\psi(t):=(\phi^m(t))'=m\phi^{m-1}(t)\phi'(t)>0$
except for some possible isolated points $\hat t$ such that
$\phi_c'(\hat t)=0$ and then $\psi(\hat t)=0$.
Using $\phi$ as a variable, we write
$\tilde\psi(\phi)=\psi(\tilde t(\phi))$ and we have
\begin{equation} \label{eq-ztildepsi-chi}
\begin{cases}
\displaystyle
\frac{d\tilde\psi}{d\phi}
=c+\chi(\hat\eta'_t(\tilde t(\phi))\phi)'_\phi
+\frac{m\phi^{m-1}(b\phi^2-a\phi)}{\tilde\psi},\\
\tilde\psi(0)=0, \quad \tilde\psi(\phi)\ge0, ~\phi\in(0,\frac{a}{b}).
\end{cases}
\end{equation}
In case (ii), we have $\tilde\psi(\frac{a}{b})=0$
and $\tilde\psi(\phi)>0$ for some left neighborhood $(\phi_1,a/b)$ of $a/b$.
Since \eqref{eq-ztilde-chi} is non-degenerate away from $\phi=0$,
the local asymptotic analysis show that there exists a $\kappa>0$ such that
$\tilde\psi(\phi)\sim \kappa(\frac{a}{b}-\phi)$, as $\phi\to (a/b)^-$.
Therefore,
$$
t^*-t_1=\int_{\phi_1}^{a/b}\frac{m\phi^{m-1}}{\tilde\psi(\phi)}d\phi
\sim
\int_{\phi_1}^{a/b}\frac{m\phi^{m-1}}{\kappa(\frac{a}{b}-\phi)}d\phi=+\infty,
$$
where $t_1=\phi_c^{-1}(\phi_1)$.
This contradicts  the case (ii).
In case (iii), we have $\tilde\psi(\phi_c(t_0))=0$
and then $\tilde\psi'_\phi(\phi_c(t_0))=0$ since $\tilde\psi(\phi)\ge0$.
There exists a constant $\omega>0$ such that
$\tilde\psi(\phi)\le \omega(\phi_c(t_0)-\phi)$ for $\phi\in(0,\phi_c(t_0))$ and
$$
t_0-0=\int_0^{\phi_c(t_0)}\frac{m\phi^{m-1}}{\tilde\psi(\phi)}d\phi
\ge \int_{\phi_c(t_0)/2}^{\phi_c(t_0)}\frac{m\phi^{m-1}}{\tilde\psi(\phi)}d\phi
\gtrsim
\int_{\phi_c(t_0)/2}^{\phi_c(t_0)}\frac{m\phi^{m-1}}{\omega(\phi_c(t_0)-\phi)}d\phi=+\infty,
$$
which contradicts  the case (iii).
Therefore, only case (i) can happen and there is no isolated points
$\hat t$ such that $\tilde\psi(\phi_c(\hat t))=0$.
That is, $t^*=+\infty$, $\tilde\psi(\phi)>0$ for all $\phi\in(0,a/b)$
and $\tilde\psi(a/b)=0$.

Let $g(\phi)\in C^1([0,a/b])$ be a nonnegative function
such that $g'(\phi)<0$ and $g(a/b)=0$.
The admissible set of functions $g$ is denoted by $D$ for convenience.
Multiplying \eqref{eq-ztildepsi-chi} by $g(\phi)$ and
integrating the resulting equation over $(0,a/b)$, we obtain
\begin{align*}
c\int_0^{a/b}g(\phi)d\phi
=&\int_0^{a/b}g(\phi)\frac{d\tilde\psi}{d\phi}d\phi
+\int_0^{a/b}g(\phi)\frac{m\phi^{m-1}(a\phi-b\phi^2)}{\tilde\psi(\phi)}d\phi
\\
&-\chi\int_0^{a/b}g(\phi)(\hat\eta'_t(\tilde t(\phi))\phi)'_\phi d\phi
\\
=&\int_0^{a/b}\Big(-g'(\phi)\tilde\psi(\phi)
+g(\phi)\frac{m\phi^{m-1}(a\phi-b\phi^2)}{\tilde\psi(\phi)}\Big)d\phi
\\
&-\chi\int_0^{a/b}-g'(\phi)\hat\eta'_t(\tilde t(\phi))\phi d\phi,
\end{align*}
since
$$
\big[g(\phi)\tilde\psi(\phi)\big]_0^{a/b}=0,
\quad
-\big[\chi g(\phi)\hat\eta'_t(\tilde t(\phi))\phi\big]_0^{a/b}=0.
$$
Therefore,
\begin{align} \nonumber
c\int_0^{a/b}g(\phi)d\phi
=&\int_0^{a/b}\Big(-g'(\phi)\tilde\psi(\phi)
+g(\phi)\frac{m\phi^{m-1}(a\phi-b\phi^2)}{\tilde\psi(\phi)}\Big)d\phi
\\ \label{eq-zcstar}
&-\chi\int_0^{a/b}-g'(\phi)\hat\eta'_t(\tilde t(\phi))\phi d\phi
\\ \nonumber
\ge&2\int_0^{a/b}
\sqrt{-g'(\phi)g(\phi)m\phi^{m-1}(a\phi-b\phi^2)}d\phi
\\ \label{eq-zvar}
&-\chi\int_0^{a/b}-g'(\phi)\hat\eta'_t(\tilde t(\phi))\phi d\phi,
\quad \forall g\in D.
\end{align}
The above equality of \eqref{eq-zvar} is attained for
$g=\hat g$ such that
\begin{equation} \label{eq-zg}
-\hat g'(\phi)\tilde\psi(\phi)
=\hat g(\phi)\frac{m\phi^{m-1}(a\phi-b\phi^2)}{\tilde\psi(\phi)}.
\end{equation}
Here we note that
$\tilde\psi(\phi)\sim \kappa(\frac{a}{b}-\phi)$,
and \eqref{eq-zg} is solvable with $\hat g(a/b)=0$
and the asymptotic behavior near zero.
Lemma \ref{le-bounds-psi} entails that $\hat g(0)<+\infty$
and $\hat g\in D$.
Hence \eqref{eq-zcstar} implies that
\begin{align*}
c=\frac{2\int_0^{a/b}
\sqrt{-\hat g'(\phi)\hat g(\phi)m\phi^{m-1}(a\phi-b\phi^2)}d\phi}
{\int_0^{a/b}\hat g(\phi)d\phi}
-\chi\frac{\int_0^{a/b}-\hat g'(\phi)\hat\eta'_t(\tilde t(\phi))\phi d\phi}
{\int_0^{a/b}\hat g(\phi)d\phi}.
\end{align*}
Now we conclude that
\begin{equation} \label{eq-zcupper}
c<\sup_{g\in D}\frac{2\int_0^{a/b}
\sqrt{-g'(\phi)g(\phi)m\phi^{m-1}(a\phi-b\phi^2)}d\phi}{\int_0^{a/b}g(\phi)d\phi}
=c_*(m,0),
\end{equation}
and according to \eqref{eq-zvar}
\begin{equation} \label{eq-zclower}
c\ge c_*(m,0)
-\chi\frac{\int_0^{a/b}-\tilde g'(\phi)\hat\eta'_t(\tilde t(\phi))\phi d\phi}
{\int_0^{a/b}\tilde g(\phi)d\phi},
\end{equation}
where the last equality of \eqref{eq-zcupper}
follows from Lemma \ref{le-var}
and $\tilde g\in D$ is the function such that
the functional in \eqref{eq-zcupper} attains $c_*(m,0)$,
whose existence is shown in Huang et al. \cite{HJMY}.
In fact, for the case without chemotaxis (i.e. $\chi=0$),
one can show that $\hat g= \tilde g$ is the
extreme function.
We note that $\tilde g$ only depends on $a,b$, and $m$,
and $\tilde g(0)$ is a fixed constant
if we normalize $\tilde g(\phi)$ such that $\int_0^{a/b}\tilde g(\phi)d\phi=1$.
Now, \eqref{eq-zclower} tells us (noticing that $\tilde g'(s)<0$ for $s\in(0,a/b)$ and
$\tilde g(a/b)=0$)
\begin{align*}
c\ge& c_*(m,0)
-\chi\int_0^{a/b}-\tilde g'(\phi)\hat\eta'_t(\tilde t(\phi))\phi d\phi
\\
\ge& c_*(m,0)
-\chi\int_0^{a/b}|-\tilde g'(\phi)|\cdot|\hat\eta'_t(\tilde t(\phi))|\phi d\phi
\\
\ge& c_*(m,0)
-\chi\int_0^{a/b}-\tilde g'(\phi)\frac{a}{2b}\frac{a}{b}d\phi
\\
=& c_*(m,0)-\chi(\tilde g(0)-\tilde g(a/b))\frac{a^2}{2b^2}
\ge \frac{2}{3}c_*(m,0),
\end{align*}
provided that $\chi\tilde g(0)\frac{a^2}{2b^2}\le \frac{1}{3}c_*(m,0)$.
$\hfill\Box$

\bigbreak

{\it \bfseries Proof of Theorem \ref{th-main}.}
We construct a mapping for finding the semi-finite TW of \eqref{eq-u2}
by solving \eqref{eq-aux-m} as following.
Define
$\overline \phi(t)=\min\{K_1t_+^\frac{1}{m-1}+C_2t_+^\frac{m}{m-1},\frac{a}{b}\}$,
and
$\underline \phi(t)=\max\{K_1t_+^\frac{1}{m-1}-C_1t_+^\frac{m}{m-1},0\}$
with $C_1$ and $C_2$ being the constants in Lemma \ref{le-bounds}.
For any $\hat\phi(t)\in \Phi$ with $\Phi$ being defined by \eqref{eq-Phi},
let $\hat\eta(t)=B[\hat\phi](t)$ as in \eqref{eq-eta},
$c=c^*(m,\chi,\hat\phi)$ as in Lemma \ref{le-cstar},
$\lambda(t)=c+\chi\hat\eta'(t)$, $\mu(t)=a-\chi\hat\eta''(t)$
and $\phi_c(t)$ be the local semi-finite type solution of \eqref{eq-aux-m}
solved in Lemma \ref{le-solv}.

Under the assumptions of Theorem \ref{th-main},
we see that
$$c^*(m,\chi,\hat\phi)\in\big[\frac{2}{3}c_*(m,0),c_*(m,0)\big)$$
and $t^*=+\infty$ according to Lemma \ref{le-varchemo},
and $\frac{1}{3}c_*(m,0)\le\lambda(t)\le\frac{4}{3}c_*(m,0)$,
$|\mu(t)|\le \frac{a}{2}$.
It follows that \eqref{eq-aux-m} is solved on $\mathbb R$
for this special $c^*(m,\chi,\hat\phi)$,
and that \eqref{eq-aux-g} and \eqref{eq-aux} are solved
for $c=c^*(m,\chi,\hat\phi)$ with the solution being exactly $\phi_c(t)$.
Denote
$$
T(\hat\phi)=\phi_c, \quad \text{~for any~} \hat\phi\in\Phi,
$$
then $T$ is well-defined and $\phi_c=T(\hat\phi)\in \Phi$
according to Lemma \ref{le-bounds}.
Note that $\Phi$ is a bounded convex subset of
the linear space $\mathscr E=C_\text{unif}^b(\mathbb R)$
endowed with the norm
$$
\|\phi\|_*=\sum_{n=1}^\infty \frac{1}{2^n}\|\phi\|_{L^\infty([-n,n])}.
$$
According to Lemma \ref{le-solv} and the energy method, we see that the solutions
are locally uniformly bounded in $C^{\gamma/m}(A)$ for some $\gamma\in(0,\frac{1}{2})$
and any compact subset $A$,
which means that $T(\Phi)$ is compact in $\mathscr E$.
The continuity of the operator $T$ follows similarly.
Applying Schauder's fixed point theorem on $\Phi$, we find that
there exists a fixed point $\phi\in\Phi$ such that $\phi=T(\phi)$.
Therefore, $\phi(t)$ is a monotonically increasing semi-finite type TW of \eqref{eq-u2}
and its speed $c=c_*(m,\chi):=c^*(m,\chi,\phi)\in\big[\frac{2}{3}c_*(m,0),c_*(m,0)\big)$
as proved in Lemma \ref{le-varchemo}.
$\hfill\Box$
\bigbreak

{\it \bfseries Proof of Theorem \ref{th-sharp}.}
Let $\phi(t)$ be a semi-finite TW of \eqref{eq-u2} and $\eta(t)=B[\phi](t)$.
According to Lemma \ref{le-solv}, we have
$$
\phi(t)=K_1t_+^\frac{1}{m-1}+K_2t_+^\frac{m}{m-1}+o(t_+^\frac{m}{m-1}),
\quad t\to 0,
$$
where $t_+=\max\{t,0\}$, and
$$
K_1=\Big(\frac{m-1}{m}(c+\chi\eta'(0))\Big)^\frac{1}{m-1},
\quad
K_2=-\frac{a-\frac{m}{m-1}\chi\eta''(0)}{2m(c+\chi\eta'(0))}K_1.
$$
It follows that $\phi\in C^1(\mathbb R)$ for $1<m<2$
and $\phi''\not\in L_\text{loc}^1(\mathbb R)$ for $m\ge2$.
The proof is completed.
$\hfill\Box$

\section*{Acknowledgement}
The research of S. Ji was
supported by NSFC Grant No. 11701184 and CSC Grant No. 201906155021.
The research of Z. Wang is supported by Hong Kong RGC GRF grant PolyU 153031/17P
(Project ID P0005368).
The research of J. Yin was supported in
part by NSFC Grant No. 11771156 and NSF of Guangzhou Grant No. 201804010391.

\end{document}